\documentclass[11pt,leqno,amssymb,a4paper]{amsart}

\usepackage{amssymb, amsmath}

\newcommand{\R}{\mathbb{R}}

\newcommand{\Z}{\mathbb{Z}}
\newcommand{\C}{\mathbb{C}}
\newcommand{\IP}{\mathbb{P}}

\newcommand{\N}{\mathbb{N}}

\newcommand{\IS}{\mathbb{S}}

\newcommand{\rs}{\mbox{$\widehat{\C}$}}

\def\SSS{{\mathcal S}}

\def\GGG{{\mathcal G}}

\def\MMM{{\mathcal M}}

\def\MMM{{\mathcal M}}

\def\YYY{{\mathcal Y}}

\def\RRR{{\mathcal R}}
\def\SSS{{\mathcal S}}
\def\UUU{{\mathcal U}}
\def\VVV{{\mathcal V}}

\newcommand{\wt}[1]{\widetilde{#1}}

\newtheorem{thm}{Theorem}[section]
\newtheorem{defn}[thm]{Definition}

\newtheorem{prop}[thm]{Proposition}

\newtheorem{lemma}[thm]{Lemma}
\newtheorem{cor}[thm]{Corollary}

\renewcommand{\qed}{\nopagebreak \begin{flushright}%end-of-proof
        \rule{2mm}{2.5mm} \end{flushright}}

\def\implies{\Rightarrow}
\newcommand{\bdry}{\partial}                     %boundary
\newcommand{\cl}{\overline}                      %closure
\def\mod{\mbox{\rm mod}}    %mod, short for modulus
\newcommand{\intersect}{\cap}                    %intersection
\newcommand{\union}{\cup}                        %union
\newcommand{\diam}{\mbox{\rm diam}}					         %diameter

\newcommand{\interior}{\mbox{int}}     %interior
\newcommand{\mtwo}[4]                            %2x2 matrices--
{\mbox{$\left(\begin{array}{cc}                  %takes four arguments
#1 & #2 \\
#3 & #4 
\end{array}
\right)$}}

\newcommand{\dettwo}[4]                          %2x2 matrices--
{\mbox{$\left|\begin{array}{cc}                  %takes four arguments
#1 & #2 \\
#3 & #4 
\end{array}
\right|$}}

\newcommand{\pf}{\noindent {\bf Proof: }}

\newcommand{\be}{\begin{enumerate}}
\newcommand{\eb}{\end{enumerate}}
\newcommand{\bi}{\begin{itemize}}
\newcommand{\ib}{\end{itemize}}
\newcommand{\bl}{\begin{list}}
\newcommand{\lb}{\end{list}}

\newcommand{\gap}{\vspace{5pt}}                 %make a space of a blank line

\newcommand{\roundness}{\mbox{\rm Round}}

\newcommand{\wtU}{{\widetilde{U}}}

\newcommand{\wtgamma}{\widetilde{\gamma}}

\newcommand{\tx}{\tilde{x}}

\newcommand{\confdim}{\mbox{\rm confdim}}
\newcommand{\hdim}{\mbox{\rm H.dim}}

\begin{document}

\author{Peter Ha\"{\i}ssinsky}
\address{LATP/CMI\\ Universit\'e de Provence\\ 39, rue Fr\'ed\'eric Joliot-Curie\\13453 Marseille cedex 13\\France}                        
\email{phaissin@cmi.univ-mrs.fr}
\author{Kevin M. Pilgrim}
\address{Dept. of Mathematics\\ Indiana University\\Bloomington\\ IN 47405\\  USA}
\email{pilgrim@indiana.edu}

\title{Thurston obstructions and Ahlfors regular conformal dimension}

%\date{\today}

\maketitle

%\vspace*{-5cm} \noindent{\it Revised version}
%\vspace*{5cm}

\begin{abstract}Let $f: S^2 \to S^2$ be a postcritically finite expanding  branched covering map of the sphere to itself.  Associated to $f$ is a canonical quasisymmetry class $\GGG(f)$ of Ahlfors regular metrics on the sphere in which the dynamics is (non-classically) conformal.   We find a lower bound on the Hausdorff dimension of metrics in $\GGG(f)$ in terms of the combinatorics of $f$. 

\gap

Soit $f:S^2\to S^2$ un rev\^etement ramifi\'e de la sph\`ere topologiquement expansif et \`a ensemble postcritique
fini. On lui associe une famille de m\'etriques Ahlfors-r\'eguli\`eres canoniques $\GGG(f)$ qui rendent
$f$ grossi\`erement conforme. On \'etablit une minoration de la dimension de Hausdorff de ces m\'etriques
en termes combinatoires. \end{abstract} 

\newpage

\section{Introduction}

A fundamental principle of dynamical systems is that in the presence of sufficient expansion, topology determines a preferred class $\GGG$ of geometric structures.    For example, suppose  $G \rightsquigarrow X$ is an action of a group $G$ on a perfect metrizable compactum $X$ by homeomorphisms.  Bowditch  \cite{bowditch:characterization} showed that if the induced diagonal action on the space of ordered triples of pairwise distinct points of $X$ is properly discontinuous and cocompact, then $G$ is hyperbolic, and there is a $G$-equivariant homeomorphism $\phi$ of $X$ onto $\bdry G$.  The boundary $\bdry G$ carries a preferred (quasisymmetry) class  of metrics in which the group elements act by uniformly quasi-M\"obius maps.  Elements of this class of metrics can be transported via $\phi$ to $X$, yielding  a class of metrics $\GGG(G \rightsquigarrow X)$ canonically associated to the dynamics in which the elements act in a geometrically special way.  

{\em Cannon's Conjecture} is equivalent to the assertion that under the hypotheses of Bowditch's theorem, whenever $X$ is homeomorphic to the two-sphere $S^2$, then the standard Euclidean metric belongs to $\GGG(G \rightsquigarrow X)$ \cite{bonk:kleiner:qsparam}.  Thus, conjecturally, no ``exotic'' metrics on the sphere arise from such group actions.  In contrast, the dynamics of certain iterated maps $f: S^2 \to S^2$ provide a rich source of examples of metrics on the sphere in which the dynamics is (non-classically) conformal.   

We now explain this precisely.  The results summarized below are consequences of the general theory developed in \cite{kmp:ph:cxci}.  
\gap

\noindent{\bf Topologically coarse expanding conformal (cxc) dynamics.}  
\begin{defn}[Topologically cxc]
A continuous, orientation-preserving, branched covering $f: S^2 \to S^2$ is called {\em topologically cxc} provided there exists a finite open covering $\UUU_0$ of $S^2$ by connected sets satisfying the following properties: 
\be  
\item[] {\bf [Expansion]}  The mesh of the covering $\UUU_n$ tends to zero as $n \to \infty$, 
where $\UUU_n$ denotes the set of connected components of $f^{-n}(U)$ as $U$ ranges over $\UUU_0$.
That is, for any finite open cover $\YYY$ of $S^2$ by open sets, there exists $N$ such that for all $n \geq N$ and all $U \in \UUU_n$, there exists $Y \in \YYY$ with $U \subset Y$.  

\item[] {\bf [Irreducibility]}  The map $f$ is {\em locally eventually onto}:  for any $x\in S^2$, and any neighborhood $W$ of $x$, there is  some $n$ with $f^n(W)=S^2$.

\item[] {\bf [Degree]} The set of degrees of maps of the form $f^k|\wtU: \wtU \to U$, where $U \in \UUU_n$, $ \wtU \in \UUU_{n+k}$, and $n$ and $k$ are arbitrary, has a finite maximum.  
\eb
We denote by $\mathbf{U}=\union_{n \geq 0}\UUU_n$.  
\end{defn}
\gap
 Note that the definition prohibits periodic or recurrent branch points, i.e.branch points $x$ for which the orbit $x, f(x), f(f(x)), \ldots$ contains or accumulates on $x$. 
 
Let $\rs$ denote the Riemann sphere.  A rational map $f: \rs \to \rs$ is a dynamical system which is conformal in the Riemannian sense.  It is called {\em semihyperbolic} if it has no parabolic cycles and no recurrent critical points in its Julia set.  A rational map $f$  which is chaotic on all of $\rs$ (that is, has Julia set the whole sphere) is  topologically cxc if and only if it is  semihyperbolic \cite[Corollary 4.4.2]{kmp:ph:cxci}.
\gap

\noindent{\bf Metric cxc.}  
Now suppose $S^2$ is equipped with a metric $d$ (that is, a distance function) compatible with its topology. 
\gap

We first recall the notion of {\em roundness}.

\gap 

\noindent{\bf Roundness.} Let $(X,d)$ be a metric space.  We denote by $B(x,r)$ and $\cl{B}(x,r)$ respectively the open and closed ball of radius $r$ about $x$.   Let $A$ be a bounded, proper  subset of $X$ with non-empty interior.  Given $x \in \interior(A)$, let 
\[ L(A,x)=\sup\{d(x,b): b \in A\}\]
and
\[ l(A,x)=\sup\{ r : r \leq L(A,x) \; \mbox{ and } \; B(x,r) \subset A\}\]
denote, respectively, the {\em outradius} and {\em inradius} of $A$ about $x$.  
While the outradius is intrinsic, the inradius depends on how $A$ sits in $X$.  The condition $r \leq L(A,x)$ is necessary to guarantee that the outradius is at least the inradius. The {\em roundness of $A$ about $x$} is defined as 
\[ \roundness(A,x) = L(A,x)/l(A,x) \in [1, \infty).\]
A set $A$ is {\em $K$-almost-round} if $\roundness(A,x) \leq K$ holds for some $x \in A$, and this 
implies that for some $s>0$, there exists a ball $B(x,s)$ satisfying  
\[ B(x,s) \subset A \subset \cl{B}(x,Ks).\]
 
\gap

\begin{defn}[Metric cxc]
A continuous, orientation-preserving branched covering $f: (S^2,d) \to (S^2,d) $ is called {\em metric cxc} provided it is topologically cxc with respect to some covering $\UUU_0$ such that there exist
\bi
\item continuous, increasing embeddings $\rho_{\pm}:[1,\infty) \to [1,\infty)$, the {\em forward and backward
roundness distortion functions}, and

\item increasing homeomorphisms $\delta_{\pm}:[0,1] \to [0,1]$, the {\em forward and backward relative
diameter distortion functions}
\ib
satisfying the following axioms:
\be
\setcounter{enumi}{3}

\item[] {\bf [Roundness distortion]} For all $n, k \in \N$ and for all
\[ U \in \UUU_n, \;\;\wtU \in \UUU_{n+k}, \;\; \tilde{x} \in \wtU, \;\;  
x
\in U\]
if
\[ f^{ k}(\wtU) = U, \;\;f^{ k}(\tilde{x}) = x \]
then the {\em backward roundness bound}
\begin{equation}
\label{eqn:backward_roundness_bound}
  \roundness(\wtU, \tilde{x}) \leq
\rho_-(\roundness(U,x))
\end{equation}
and the {\em forward roundness bound}
\begin{equation}
\label{eqn:forward_roundness_bound}
\roundness(U, x) \leq
\rho_+(\roundness(\wtU,
\tilde{x})).
\end{equation}
hold.  

In other words: for a given element of $\mathbf{U}$,
iterates of $f$ both forward and backward distorts its roundness by an
amount independent of the iterate.

\item[] {\bf [Diameter distortion]} For all $n_0, n_1, k \in \N$ and for all
\[ U \in \UUU_{n_0}, \;\;U' \in \UUU_{n_1}, \;\;\wtU \in \UUU_{n_0+k},
\;\;\wtU'
\in
\UUU_{n_1+k}, \;\; \wtU' \subset \wtU, \;\; U' \subset U\]
if
\[ f^k(\wtU) = U, \;\;f^k(\wtU') = U'\]
then
\[ \frac{\diam\wtU'}{\diam\wtU} \leq \delta_-\left(\frac{\diam U'}{\diam
U}\right)\]
and
\[ \frac{\diam U'}{\diam U} \leq \delta_+\left(\frac{\diam \wtU'}{\diam
\wtU}\right)\]
hold.

In other words:  given two nested elements of $\mathbf{U}$, iterates of  
$f$
both forward and backward distort their relative sizes by an amount
independent of the iterate.
\eb
\end{defn}
\gap

 A homeomorphism $h: X \to Y$ between metric spaces is called {\em quasisymmetric} provided 
 there exists a homeomorphism $\eta: [0,\infty) \to [0,\infty)$ such that $d_X(x,a) \leq td_X(x,b)$
 implies $d_Y(f(x),f(a)) \leq \eta(t)d_Y(f(x),f(b))$ for all triples of points $x, a, b \in X$ and all $t \geq 0$.   Loosely:  $h$ distorts ratios of distances, and the roundness of balls, by controlled amounts.  
 
An orientation-preserving  branched covering map $f: \IS^2 \to \IS^2$  from the standard 
Euclidean sphere to itself is metric cxc if and only if it is quasisymmetrically conjugate 
to a semihyperbolic rational map with Julia set the whole sphere 
\cite[Theorems 4.2.4 and 4.2.7]{kmp:ph:cxci}.   
The class of metric cxc dynamical systems is closed under quasisymmetric conjugation, 
and a topological conjugacy between metric cxc maps is quasisymmetric  
\cite[Theorem 2.8.2]{kmp:ph:cxci}.

\gap

\noindent{\bf Conformal gauges.} The {\em conformal gauge} of a metric space $X$ is the set of all metric spaces quasisymmetrically equivalent to $X$.  A metric space $X$ is {\em Ahlfors regular  of dimension
$Q$} provided there is a Radon measure $\mu$ and a constant $C\geq 1$ such that
 for any $x\in X$ and $r\in (0,\diam X]$,
$$\frac{1}{C}r^Q \leq \mu(B(x,r)) \leq C r^Q\,.$$
The Hausdorff dimension $\hdim(X)$ of an Ahlfors $Q$-regular metric space $X$ is equal to $Q$.
  
\gap

Suppose now that $f: S^2 \to S^2$ is topologically cxc.  By \cite[Corollary 3.5.4]{kmp:ph:cxci} we have 

\begin{thm}[Canonical gauge]
Given a topologically cxc dynamical system $f: S^2 \to S^2$, there exists an Ahlfors regular metric $d$ on $S^2$, unique up to quasisymmetry, such that $f:(S^2,d) \to (S^2,d)$ is metrically cxc.   
\end{thm}

It follows that the set $\GGG(f)$ of all Ahlfors regular metric spaces $Y$ quasisymmetrically equivalent to $(S^2,d)$ is an invariant, called the {\em Ahlfors regular conformal gauge},  of the topological conjugacy class of $f$.  
Therefore, the {\em Ahlfors regular conformal dimension}
\[ \confdim_{AR}(f) := \inf_{Y \in \GGG(f)} \hdim(Y)\]
is a numerical topological dynamical invariant as well; it is distinct from the entropy.  Moreover, this invariant almost characterizes  rational maps among topologically cxc maps on the sphere.
In \cite[Theorem 4.2.11]{kmp:ph:cxci} the following theorem is proved.

\begin{thm}[Characterization of rational maps]
\label{thm:char_of_rm}
A topologically cxc map $f: S^2 \to S^2$  is topologically conjugate to a semihyperbolic rational map if and only if the Ahlfors regular conformal dimension $\confdim_{AR}(f)$ is equal to $2$, and is achieved by an Ahlfors regular metric.
\end{thm}

There are many examples of topologically cxc maps which are not topologically conjugate to rational maps.   The following are well-known combinatorial obstructions.  
\gap

\noindent{\bf Thurston obstructions.} Let $f: S^2 \to S^2$ be an orientation-preserving branched covering.  The Riemann-Hurwitz formula implies that the cardinality of the set $C_f$ of branch points at which $f$ fails to be locally injective is equal to $2\deg(f)-2$, counted with multiplicity, where $\deg(f)$ is the degree of $f$.  The {\em postcritical set} is defined as $P_f = \cl{\union_{n>0}f^{n} (C_f)}$.  Under the assumption that the postcritical set is finite, Thurston characterized when  $f$ is  equivalent to a rational map $R$ in the following sense:  $h_0\circ f = R \circ h_1$ for orientation-preserving homeomorphisms $h_0, h_1$ which are homotopic through homeomorphisms fixing $P_f$ pointwise; see \cite{DH1}.  The obstructions which arise are of the following nature.  

A {\em multicurve} $\Gamma \subset S^2-P_f$ is a finite set of simple, closed, unoriented curves 
\[ \Gamma=\{\gamma_1, \gamma_2, \ldots, \gamma_m\}\]
in $S^2-P_f$ satisfying the following properties:  (i) they are disjoint and pairwise distinct, up to free homotopy in $S^2-P_f$, and (ii) each curve $\gamma_j$ is {\em non-peripheral}---that is, each component of $S^2-\gamma_j$ contains at least two elements of $P_f$.   
A multicurve $\Gamma$ is {\em invariant} if for each $\gamma_j \in \Gamma$, every connected component $\delta$ of $f^{-1}(\gamma_j)$ is either homotopic in $S^2-P_f$  to an element $\gamma_i \in \Gamma$, or else is peripheral.  If $\alpha$ and $\beta$ are unoriented curves in $S^2-P_f$, we write $\alpha \sim \beta$ if they are freely homotopic in $S^2-P_f$. 

Let $\Gamma$ be an arbitrary multicurve  and $Q\ge 1$. Let $\R^{\Gamma}$ denote 
the real vector space with basis $\Gamma$; thus $\gamma_j$ is identified with the $j$th standard basis vector of $\R^{\# \Gamma}$. Define 
\[ f_{\Gamma, Q}: \R^\Gamma \to \R^\Gamma\]
by 
\[ f_{\Gamma, Q}(\gamma_j) =\sum_{\gamma_i\in\Gamma} \sum_{\delta \sim \gamma_i} |\deg(f:\delta \to \gamma_j)|^{1-Q}\gamma_i.\]
In words:  the $(i,j)$-matrix coefficient $(f_{\Gamma,Q})_{i,j}$ is obtained by considering 
the connected preimages $\delta$ of $f^{-1}(\gamma_j)$ homotopic to $\gamma_i$ in $S^2-P_f$; 
 recording the positive degree of the restriction $f|_\delta: \delta \to \gamma$, raised to the power $(1-Q)$; and summing these numbers together.  If there are no such curves $\delta$, the coefficient is defined to be zero. Note that in the definition, we do not require invariance.  
 
Since $f_{\Gamma, Q}$ is represented by a non-negative matrix,  one can apply the structure theory for such matrices, summarized at the beginning of Appendix \ref{secn:appendix}.  This theory implies the following results.  

The matrix $(f_{\Gamma,Q})$ has a real non-negative Perron-Frobenius eigenvalue $\lambda(f_{\Gamma,Q})$ equal to its spectral radius and a corresponding non-negative  eigenvector $v(f_{\Gamma,Q})$. A multicurve $\Gamma$ is called {\em irreducible} if given any $\gamma_i, \gamma_j \in \Gamma$ there exists an iterate $q\geq 1$ such that the corresponding coefficient $(f_{\Gamma,Q}^q)_{i,j}$ is positive; this property is independent of $Q$.  For an irreducible multicurve, the Perron-Frobenius eigenvalue  is positive, has geometric multiplicity one, and is strictly larger than the norm of all other eigenvalues; the corresponding eigenvector is also strictly positive.   

A {\em Thurston obstruction} is defined as a multicurve $\Gamma=\{\gamma_1, \ldots, \gamma_m\}$ 
for which the inequality $\lambda(f_{\Gamma,2}) \geq 1$ holds.   An obstruction always contains an 
irreducible obstruction with the same Perron-Frobenius eigenvalue (cf. Appendix \ref{secn:appendix}).  
By a theorem of McMullen \cite{ctm:renorm}, a semihyperbolic rational map has no obstructions unless it is extremely special (see below). 
 
The reason these form obstructions to (classical Riemannian) conformality is roughly the following; 
see \cite{DH1} for details. Suppose a semihyperbolic rational map $f: \rs \to \rs$ had an obstruction.  
Then one could find a collection of disjoint annular neighborhoods $A_j$ of $\gamma_j$ such that the 
vector of classical moduli $(\mod(A_1), \ldots, \mod(A_m))$ is a scalar multiple of a Perron-Frobenius eigenvector $v$.  
Classical moduli are subadditive and monotone:  the sharp Gr\"otzsch inequality implies that 
if $A(\delta_k), k=1, \ldots, l$ are disjoint essential open subannuli of $A_i$, then 
$\sum_{k=1}^l \mod(A(\delta_k)) \leq \mod(A_i)$; equality holds if and only if each $A(\delta_k)$ is a right  
Euclidean subannulus in a conformally equivalent Euclidean metric on $A_i$,  
and the union of their closures contains $A_i$.   
If $f: A(\delta_k) \to A_j$ is a degree $d$ covering, then $\mod(A(\delta_k))=\mod(A_j)/d$.   
It follows by induction that for fixed $j \in \{1, \ldots, m\}$ and for all $n \in \N$,  
the $j$-th coordinate of the vector $f^n_{\Gamma,2}(v)$ is a lower bound for the maximum modulus 
of an annulus homotopic in $\IP^1-P_f$ to $A_j$.  
It follows that such a rational map cannot have an obstruction unless 
it is extremely special---a so-called {\em integral Latt\`es example} \cite{DH1}.    
In this case, $\#P_f=4$, $(f_{\Gamma, 2})=(1)$, and $f$ lifts under a twofold covering ramified at 
$P_f$ to an unbranched covering map of the complex torus given by $z \mapsto dz$ in the group law, where $d=\deg(f)$.   
Summarizing, we say that branched covering $f$ is {\em obstructed} if (i) 
it is not topologically conjugate to an integral Latt\`es example, and (ii) it has an obstruction.

Suppose that the map $f: S^2 \to S^2$ is topologically cxc.  Since the property of being obstructed is invariant under topological conjugacy, 
Theorem \ref{thm:char_of_rm} yields
\begin{center}
$f$ obstructed $\implies \confdim_{AR}(f)$ is either 
$\left\{\begin{array}{rl}
2, & \mbox{but not realized, or}\\
>   2. & \; 
\end{array}\right.
$
\end{center}

Our main result, which was inspired by discussions with M. Bonk and L. Geyer, quantifies the influence of  obstructions  on the Ahlfors regular conformal gauge $\GGG(f)$.  

Let $\Gamma$ be a multicurve.  If $\Gamma$ contains an irreducible multicurve, then there is a unique value $Q(\Gamma)\geq 1$ such that $\lambda(f_{\Gamma,Q(\Gamma)})=1$ (Lemma \ref{lemma:decreasing}). Otherwise, we set $Q(\Gamma)=0$.
Define 
\[ Q(f) =\sup\{Q(\Gamma): \Gamma \  \mbox{is a multicurve} \}.\]

We note that:
\begin{itemize}
\item If $f$ is obstructed, then $Q(f) \geq 2$ and is a rough numerical measurement of the extent to which $f$ is obstructed.   
\item If $\#P_f<\infty$, then there are only finitely many possible irreducible matrices $f_{\Gamma,2}$, and the supremum is achieved by some multicurve.  
\item If $f$ is not conjugate to a Latt\`es example, and if 
$\Gamma$ is not an obstruction, then $Q(\Gamma)<2$, by Lemma \ref{lemma:decreasing}.
\end{itemize}

We prove:  

\begin{thm}  
\label{thm:ARconfdim}
Suppose $f: S^2 \to S^2$ is topologically cxc.  Then 
\[ \confdim_{AR}(f) \geq Q(f).\]
\end{thm}

The finite subdivision rules of Cannon, Floyd, and Parry \cite{cfp:fsr} provide a wealth of examples of topologically cxc maps on the sphere  \cite[\S\,4.3]{kmp:ph:cxci}.
As a special case of the above theorem, we have the following. 

\begin{cor}
Suppose $\RRR$ is a finite subdivision rule with bounded valence, mesh going to zero, underlying surface the two-sphere, and whose subdivision map $f: S^2 \to S^2$ is orientation-preserving.  Then $\confdim_{AR}(f) \geq Q(f)$.  
\end{cor}

In \cite[Conjecture 6.4]{bonk:icm:qcgeom} it is guessed that for obstructed maps  
induced by such finite subdivision rules, equality actually holds.  
The preceding corollary establishes one direction of this conjecture.  
Our methods are in spirit similar to those sketched above for the classical case $Q=2$.  
Instead of classical analytic moduli, combinatorial moduli are used.   
The outline of our argument is the same as the brief sketch in 
\cite{bonk:icm:qcgeom}.  However, Theorem \ref{thm:ARconfdim} applies to maps 
which need not be postcritically finite and hence need not arise from finite 
subdivision rules.  A key ingredient is the construction 
of a suitable metric on $S^2$ in which the coverings $\UUU_n$ have the geometric regularity property of being a family of 
{\em uniform quasipackings}.   
Also, our proof  makes use of a succinct comparison relation 
(Proposition \ref{prop:comb_mod} below) between combinatorial and analytic 
moduli articulated by the first author in \cite{ph:emp}. 

Unfortunately our proof is somewhat indirect:  apart from the bound on dimension, our methods shed very little light on the structure of the elements of the gauge $\GGG(f)$.    

%\gap

\newpage

\noindent{\bf Outline of paper.}

In \S 2 we develop the machinery of combinatorial $Q$-moduli of path families associated to sequences $(\SSS_n)_n$  of coverings of surfaces.  Much of this material is now standard.   

In \S 3 we state results that relate combinatorial and analytic moduli in Ahlfors regular metric spaces.   These results apply to covering sequences $(\SSS_n)_n$ which are {\em quasipackings with mesh tending to zero}.   

In \S 4, we briefly recall the construction in \cite{kmp:ph:cxci} of the gauge $\GGG(f)$ and its properties. We also prove that when $S^2$ is equipped with any metric in the gauge $\GGG(f)$, the sequence of coverings $(\UUU_n)_n$ defines a uniform sequence of quasipackings. 

In \S 5, we complete the proof of Theorem \ref{thm:ARconfdim}.

In Appendix A, we summarize facts about non-negative matrices and prove Lemma \ref{lemma:decreasing}.
\gap

\noindent{\bf Acknowledgments.} We thank M.\,Bourdon, A.\,Ch\'eritat and J.\,Rivera-Letelier 
for their comments on an earlier version of this manuscript.  
We also thank the anonymous referee for comments that
substantially improved the clarity of our exposition.  Both authors were partially supported
by the  project ANR ``Cannon'' (ANR-06-BLAN-0366).  The second author was supported by NSF grant DMS-0400852.
\gap

\noindent{\bf Notation.}  For positive quantities $a, b$, we write $a \lesssim b$ (resp. $a \gtrsim b$) if there is a universal constant $C>0$ such that $a \leq Cb$ (resp. $a \geq Cb$).  The notation $a \asymp b$ will mean $a \lesssim b$ and $a \gtrsim b$.  

If $A$ is a matrix, the notation $A \geq 0$ means the entries of $A$ are non-negative, and $A \geq B$ means $A-B \geq 0$. 

The cardinality of a set $A$ is denoted $\#A$.

\section{Combinatorial moduli}

\noindent{\bf Definitions.}
\gap
Let $\SSS$ be a covering of a  topological
space $X$, and let $Q\ge 1$.
Denote by $\MMM_Q(\SSS)$ the set
of functions $\rho:\SSS\to \R_+$ such that $0<\sum\rho(s)^Q<\infty$; 
elements of $\MMM_Q(\SSS)$ we call  {\sl admissible metrics}.
For $K\subset X$ we denote by $\SSS(K)$ the set of elements of $\SSS$ which intersect $K$.  The {\sl $\rho$-length} of $K$ is by definition
$$\ell_\rho
(K)=\sum_{s\in\SSS( K)} \rho(s)$$
and its {\em $\rho$-volume} is
$$V_{\rho}(K)=\sum_{s\in\SSS(K)} \rho(s)^Q\,.$$
If $\Gamma$ is a family of curves in $X$ and if $\rho\in\MMM_Q(\SSS) 
$, we define
$$L_\rho(\Gamma,\SSS)=\inf_{\gamma\in\Gamma} \ell_\rho(\gamma),$$
$$\mod_Q(\Gamma, \rho, \SSS)=\frac{V_{\rho}(X)}{L_\rho(\Gamma,\SSS)^Q},$$ 
and the {\em combinatorial modulus} by
$$\mod_Q(\Gamma,\SSS) = \inf_{\rho\in\MMM_Q(\SSS)}
\mod_Q(\Gamma,\rho,\SSS).$$

A metric $\rho$ for which $\mod_Q(\Gamma,\rho,\SSS)=\mod_Q(\Gamma, \SSS)$ will be called {\em optimal}.  We will consider here only finite coverings; in this case the proof of the existence of optimal metrics is a straightforward argument in linear algebra.  
The following result is the analog of the classical Beurling's criterion
which  characterises optimal metrics.

\begin{prop}\label{prop:metric_opt}
Let $\SSS$ be a finite cover of a space $X$,
$\Gamma$ a family of curves and $Q>1$.
An admissible metric $\rho$ is optimal if and only if there is a non-empty subfamily $ 
\Gamma_0\subset\Gamma$ and non-negative scalars $\lambda_{\gamma}$, $ 
\gamma\in\Gamma_0$, such that \be
\item for all $\gamma\in\Gamma_0$, $\ell_\rho(\gamma) = L_\rho(\Gamma, 
\SSS)$\,;
\item for any $s\in\SSS$, $$Q\rho(s)^{Q-1}=\sum \lambda_{\gamma}$$
where the sum is taken over curves in $\Gamma_0$ which  
go through $s$.\eb
Moreover, an optimal metric is unique up to scale.\end{prop}

For a proof, see Proposition 2.1 and Lemma 2.2 in \cite{ph:emp}.
\gap

\noindent{\bf Monotonicity and subadditivity.}
\gap

\begin{prop}\label{prop:subadditivity}
Let $\SSS$ be a locally finite cover of a topological space $X$ and $Q \geq 1$.

\be
\item If $\Gamma_1 \subset \Gamma_2$ then $\mod_Q(\Gamma_1,\SSS) \leq \mod_Q(\Gamma_2, \SSS)$.
\item Let $\Gamma_1,\ldots,\Gamma_n$ be a set of curve families in $X$ and $Q\ge 1$.
Then $$\mod_Q(\cup\Gamma_j,\SSS)\le \sum\mod_Q(\Gamma_j,\SSS)\,.$$
Furthermore, if $\SSS(\Gamma_i)\intersect\SSS(\Gamma_j)=\emptyset$ for $i \neq j$, then 
$$\mod_Q(\cup\Gamma_j,\SSS)=\sum\mod_Q(\Gamma_j,\SSS)\,.$$
\eb
\end{prop}

The proof is the same as the standard one for classical moduli (see for instance \cite[Thms 6.2 and 6.7]{vaisala:lectures_qc}) and so is omitted.

\gap

\noindent{\bf Naturality under coverings.}

A closed (resp. open) {\em annulus} in a surface $X$ is a subset homeomorphic to $[0,1]\times S^1$ (resp. $(0,1)\times S^1$).  
Suppose $A$ is an annulus in a surface $X$ and $\SSS$ is a finite covering of $A$ by subsets of $X$.  
For $Q \geq 1$ we define 
\[ \mod_Q(A,\SSS)=\mod_Q(\Gamma, \SSS)\]
where $\Gamma$ is the set of closed curves which are contained in $A$ and which separate the boundary components of $A$.  

Note that $\mod_Q(A, \SSS)$ is an invariant of the triple $(X, A, \SSS)$ and is not purely intrinsic to $A$. 
The following result describes how combinatorial moduli of annuli change under coverings.  
Since the elements of $\SSS$ meeting $A$ need not be contained in $A$, it is necessary to have some additional 
space surrounding $A$ on which the covering map is defined.

\begin{prop}\label{prop:covering} Suppose $A, B, A', B'$ are open annuli such that $\cl{A} \subset B$, 
$\cl{A'} \subset B'$, $A$ is essential in $B$, and $A'$ is essential in $B'$.  
Let $f:B'\to B$ be a covering map
of degree $d$ such that $f|A': A' \to A$ is also a covering map of degree $d$.
Let $\SSS$ be a finite cover of $A$ by Jordan domains $s \subset B$ and $\SSS'$ be
the induced covering of $A'$, i.e.  the covering whose elements $s'$ are the components of $f^{-1}(\{s\})$, $s\in\SSS$. 
 Then, for $Q>1$,
$$\mod_Q(A',\SSS') = d^{1-Q}\cdot \mod_Q(A,\SSS)\,.$$
\end{prop}

\pf Let $\Gamma, \Gamma'$ denote respectively the curve families in $A, A'$ separating the boundary components.  We note that since $f$ is a covering and each piece of $\SSS$ is a Jordan domain,
$f^{-1}(s)$ has $d$ components each of which is also a Jordan domain.

Let $\rho$ be an optimal metric for $\mod_Q(\Gamma,\SSS)$.
Consider the subfamily $\Gamma_0$ and the scalars $ 
\lambda_\gamma$
given by Proposition \ref{prop:metric_opt}.
Set $\Gamma_0'=f^{-1}(\Gamma_0)$, $\rho'=\rho\circ f$, and for $\gamma' \in \Gamma_0'$ define 
$\lambda_{\gamma'}= \lambda_{f(\gamma')}$.

The preimage $\gamma'$ of a curve  $\gamma$ in $\Gamma$ is connected and all the preimages 
of $s\in\SSS(\gamma)$ belong to $\SSS'(\gamma')$.
Therefore, for $\gamma' \in\Gamma_0'$, 
one has $\ell_{\rho'}(\gamma')=d L_\rho (\Gamma)$,
and for any other curve,  $\ell_{\rho'}(\gamma')\ge d L_\rho(\Gamma)$.

Clearly, for any $s'\in\SSS'$, $$Q\rho'(s')^{Q-1}=\sum_{\gamma'\in\Gamma'_0} \lambda_ 
{\gamma'}$$
so that Proposition \ref{prop:metric_opt} implies that $\rho'$ is
optimal.

It follows that
$$\mod_Q(\Gamma',\SSS') =\frac{d V_Q(\rho)}{(d L_\rho(\Gamma))^Q}=  
d^{1-Q}\cdot\mod_Q(\Gamma,\SSS)\,.$$
\qed

\section{Combinatorial moduli and Ahlfors regular conformal dimension}

Under suitable conditions, the combinatorial moduli obtained from a sequence $(\SSS_n)_n$ of coverings can be used to approximate
analytic moduli on metric measure spaces. 
Suppose $(X,d,\mu)$ is a metric measure space,
$\Gamma$ is a family of curves in $X$, and $Q\ge 1$.  The {\em (analytic) $Q$-modulus} of $\Gamma$ is defined by 
$$\mod_Q(\Gamma)=\inf \int_X \rho^Q d\mu$$
where the infimum is taken over all measurable functions
$\rho:X\to\R_+$ such that $\rho$ is {\em admissible}, i.e. $$\int_{\gamma}\rho ds\ge 1$$ for all
$\gamma\in \Gamma$ which are rectifiable. 
If $\Gamma$ contains no rectifiable curves, $\mod_Q(\Gamma)$ is defined to be zero.  
Note that when $\Gamma$ contains a constant curve, then
there are no admissible $\rho$, so we set $\mod_Q\Gamma=+\infty$.
When $X \subset \C$ is a domain, $\mu$ is Euclidean area,  
and $Q=2$, this definition coincides with the classical one.  

\gap
The approximation result we use requires the sequence of coverings $(\SSS_n)_n$ to 
be a {\em uniform family of quasipackings}.  

\gap
\begin{defn}[Quasipacking] A {\em quasipacking} of a metric space is
a locally finite cover $\SSS$ such that there is some constant
$K\ge 1$ which satisfies the following property.
For any $s\in\SSS$, there are two balls $B(x_s, r_s)\subset s\subset B(x_s, K\cdot r_s)$ 
such that the family $\{B(x_s, r_s)\}_{s\in\SSS}$ consists of pairwise disjoint balls.  
A family $(\SSS_n)_n$ of quasipackings is called {\em uniform } if the mesh of $\SSS_n$ tends to zero as $n \to \infty$ and the constant $K$ defined above can be chosen independent of $n$. 
\end{defn}

The next result says roughly that for the family consisting of all sufficiently large curves, analytic and combinatorial moduli are comparable.  

\begin{prop}\label{prop:comb_mod} Suppose $Q > 1$, $X$ is an Ahlfors $Q$-regular compact  
metric space, and $(\SSS_n)_n$ is a sequence of uniform quasipackings.  Fix $L>0$, and let
 $\Gamma_L$ be the family of curves in $X$ of diameter at least $L$.
Then either 
\be
\item $\mod_Q(\Gamma_L)=0$ and $\lim_{n \to \infty}\mod_Q(\Gamma_L,\SSS_n)=0$, or 
\item $\mod_Q(\Gamma_L)>0$, and there exists constants $C \geq 1$ independent of $L$ and $N=N(L) \in \N$ such that  for any  $n>N$,
$$\frac{1}{C}\mod_Q(\Gamma_L,\SSS_n)\leq \mod_Q(\Gamma_L)\leq C\mod_Q(\Gamma_L,\SSS_n).$$
\eb
\end{prop}

See Proposition B.2 in \cite{ph:emp}. 
\gap

\begin{cor}\label{cor:dim} 
Under the hypotheses of Proposition \ref{prop:comb_mod}, if  
$Q > \confdim_{AR}(X)\geq 1$, and if $\Gamma$ is a curve family contained in some $\Gamma_L$, $L>0$, then $\lim_{n \to \infty}\mod_Q(\Gamma,\SSS_n)=0$.
\end{cor}

In other words, if there is a curve family $\Gamma$ each of whose elements has diameter at least $L>0$, and if for some $Q>1$  
one has  $$\mod_Q(\Gamma,\SSS_n)\gtrsim 1$$
for all $n$, then the AR-conformal dimension is at least $Q$.

\gap
The lower bound on the diameter is necessary.   
If $\Gamma$ is any family of curves such that, for any $n$, there is some $\gamma\in\Gamma$ contained 
in an element of $\SSS_n$, then $\mod_Q(\Gamma,\SSS_n)\ge 1$ for all $n$.  
So, for example, 
if $\Gamma$ consists of a countable family of non-constant curves as above, 
then $\mod_Q(\Gamma)=0$ while $\mod_Q(\Gamma, \SSS_n) \ge 1$ for all $n$.

\gap

\pf By assumption, there is some  metric $d$ in the conformal gauge of $X$ which is Ahlfors regular 
of dimension $p \in (\confdim_{AR}X, Q)$, and $\Gamma \subset \Gamma_L$ for some $L>0$.  
Let $d'=d^{p/Q}$.  Then the sequence $\{\SSS_n, n\geq 0\}$ is again a family of uniform quasipackings.  
Though $d'$ is Ahlfors regular of dimension $Q$, it has no rectifiable curves.  
In particular, for the metric $d'$, we have $\mod_Q(\Gamma_L)=0$.  
By Proposition \ref{prop:comb_mod}, $\mod_Q(\Gamma_L, \SSS_n) \to 0$ as $n \to \infty$.
Since $\Gamma \subset \Gamma_L$, the monotonicity of moduli (Proposition \ref{prop:subadditivity}) 
implies that $\mod_Q(\Gamma,\SSS_n)\le \mod_Q(\Gamma_L,\SSS_n)$ so that
 $\mod_Q(\Gamma, \SSS_n)$ tends to $0$ as well.
\qed

{\noindent\bf Remark.}  Corollary \ref{cor:dim} takes its origin in the work of Pansu
\cite[Prop.\,3.2]{pansu:confdim} where a similar statement is proved for his {\it modules grossiers}.
It is also closely related  to a theorem of Bonk and Tyson 
which asserts that if the $Q$-modulus of curves in a
$Q$-Ahlfors regular space is positive, then the Ahlfors regular conformal dimension of that space
is $Q$ \cite[Thm\,15.10]{heinonen:analysis}.  In particular, if a metric is $Q$-regular for some $Q$ strictly larger than the Ahlfors regular conformal dimension,
then the $Q$-modulus of any non-trivial family of curves is zero.   

\section{The conformal gauge of a topological cxc map}

In this section, we recall from \cite{kmp:ph:cxci} the construction of the metrics associated to topologically cxc maps, 
specialized to the case of maps $f: S^2 \to S^2$.   
After summarizing their properties, we prove that with respect to these metrics, the induced coverings $\UUU_n$ obtained by 
pulling back an initial covering $\UUU_0$ under iteration form a sequence of uniform quasipackings.
\gap

\noindent{\bf Associated graph $\Sigma$.}  Suppose $f: S^2 \to S^2$ is topologically cxc with respect to an open covering $\UUU_0$.  Let $\Sigma$ be the graph whose vertices are elements of $\union_n \UUU_n$, together with a distinguished root vertex $o={S^2}=\UUU_{-1}$.  The set of edges is defined as a disjoint union of two types of edges:  horizontal edges join elements $U_1, U_2 \in \UUU_n$ if and only if $U_1 \intersect U_2 \neq \emptyset$, while vertical edges join elements $U \in \UUU_n, V\in \UUU_{n\pm 1}$ at consecutive levels if and only if $U \intersect V \neq \emptyset$.  Note that there is a natural map $F: \Sigma \to \Sigma$ which is cellular on the complement of the set of closed edges meeting $\UUU_0$.  
\gap

\noindent{\bf Associated metrics.}  Equip $\Sigma$ temporarily with the length metric
$d(\cdot, \cdot)$ in which edges are isometric to unit intervals.  

Axiom [Expansion] implies that 
the metric space $\Sigma$ is hyperbolic in the sense of Gromov \cite[Theorem 3.3.1]{kmp:ph:cxci}; 
see \cite{ghys:delaharpe:groupes} for background on hyperbolic metric spaces.  
One may define its compactification in the following way.

 Fix $\varepsilon > 0$.  For $\xi \in \Sigma$ 
let $\varrho_\varepsilon (\xi)=\exp(-\varepsilon d(o,\xi))$.  Define  
a new metric $d_\varepsilon$ on $\Sigma$ by
\[ d_\varepsilon(\xi,\zeta) = \inf \ell_\varepsilon(\gamma)\]
where 
\[ \ell_\varepsilon(\gamma) = \int_\gamma \varrho_\varepsilon \ ds\]
and where as usual the infimum is over all rectifiable curves in the metric space $(\Sigma, d)$ joining $\xi$ to $\zeta$. 
The resulting metric space $\Sigma_\varepsilon=(\Sigma, d_\varepsilon)$ is incomplete. Its complement in its
completion defines the boundary $\bdry_\varepsilon \Sigma$ 
which is an Ahlfors regular metric space of dimension $\frac{1}{\varepsilon}\log \deg(f)$ by axiom [Degree] if $\varepsilon$ is small enough.  The map $F$ is $e^\varepsilon$-Lipschitz in the $d_\varepsilon$-metric, so it extends to $\bdry_\varepsilon\Sigma$.

If $\varepsilon$ is sufficiently small, the boundary $\bdry_\varepsilon\Sigma$ is homeomorphic to the usual Gromov boundary, and there is a natural homeomorphism $\phi: S^2 \to \bdry_\varepsilon\Sigma$ given as follows.  For $x \in S^2$ let $U_n(x)$ be any element of $\UUU_n$ containing $x$.  We may regard $U_n(x)$ as a vertex of $\Sigma$ and hence as an element of its completion $\cl{\Sigma}_\varepsilon$.  The definitions of $\Sigma$ and of $d_\varepsilon$ imply that 
\[ \phi(x) = \lim_{n \to \infty} U_n(x) \in \bdry_\varepsilon\Sigma\]
exists and is independent of the choice of sequence $\{U_n(x)\}_n$.   The homeomorphism $\phi$  
conjugates $f$ on $S^2$ to the map $F$ on $\bdry_\varepsilon \Sigma$.
\gap

\noindent{\bf Associated metrics on $S^2$.}  
 {\em A priori} the boundary $\bdry_\varepsilon\Sigma$ depends on the choice of $\UUU_0$ and of $\varepsilon$.  However, by \cite[Proposition 3.3.12]{kmp:ph:cxci}, its quasisymmetry class is independent of such choices, provided the covering satisfies axiom [Expansion] and the parameter is small enough to guarantee that $\phi$ is a homeomorphism.  We remark that balls for  such metrics need not be connected.  

The {\em Ahlfors regular conformal gauge} $\GGG(f)$ is then defined as the set of all Ahlfors regular metrics on $S^2$ quasisymmetrically equivalent to a metric of the form $\phi^*(d_\varepsilon)$.  Elements of $\GGG(f)$ will be referred to as  {\em associated metrics}.  

It what follows, for convenience we denote by $d_\varepsilon$  the pulled-back metric $\phi^*(d_\varepsilon)$ on $S^2$.  

\begin{thm}\label{thm:qpack_exists}
Let $f:S^2\to S^2$ be a topological cxc dynamical system with respect to an open covering $\VVV_0$.  Then there exists a finite 
cover $\UUU_0$ of $S^2$ by Jordan domains such that, for any associated metric, the sequence of coverings $\{\UUU_n, n \geq 0\}$ is a  uniform family of quasipackings.  
\end{thm}

\pf It is easily shown that the property of being a uniform quasipacking is preserved under quasisymmetric changes of metric.  Hence, it suffices to show the conclusion for a metric $d_\varepsilon$ as constructed above.

For convenience, equip $S^2$ with the standard Euclidean spherical metric and denote the 
resulting metric space by $\IS^2$.  Then small spherical balls $D(x,r)$ are Jordan domains.  
For each $x\in \IS^2$, consider an open ball $U_x=D(x,r_x)$ centred at $x$.  
By expansion, there exists $n_0$ such that no element of $\VVV_n$, $n \geq n_0$, contains more than one critical value of $f$.   By choosing $r_x$ sufficiently small and sufficiently generic, we may arrange so that each $U_x$ (i) is contained in some element of $\VVV_{n_0}$, and (ii) does not contain a critical value of an iterate of $f$ on its boundary.  A covering of a disk ramified above at most one point is again a disk, by the Riemann-Hurwitz formula.  It follows that every iterated preimage of $U_x$ under $f$ is a Jordan domain. 
 
Let $\UUU_0=\{U_{x_j}\}_j$ be a finite subcover. From the $5r$-covering
theorem \cite[Thm 1.2]{heinonen:analysis}, we may assume that the balls 
$\{D(x_j,r_j/5)\}_j$ are pairwise disjoint.
Since we assumed that each disk in $\UUU_0$ was contained in an element of $\VVV_{n_0}$,
it follows that $\UUU_0$ satisfies both axioms [Expansion] and [Degree].
Furthermore, there is some $r_0>0$ such
that the collection $\{B_{\varepsilon}(x_j,r_0)\}$ of balls in the  metric $d_\varepsilon$ is a disjointed family.   

The proof is completed by appealing to the axiom [Expansion] and  to the fact that with respect to the metric $d_\varepsilon$, iterates of $f$ distort 
 elements of $\mathbf{U}=\union_n \UUU_n$ by controlled amounts.  
More precisely, it follows from \cite[Prop. 3.3.2]{kmp:ph:cxci} that there is a constant $C\ge 1$  for which the following property holds.
If $\tx\in\wtU\in\UUU_{n}$, $f^n(\tx)= x_j$, $f^n(\wtU)= U_j$,
then $$B_{\varepsilon}(\tx, (r_0/C)e^{-\varepsilon n})\subset
\wtU\subset B_{\varepsilon}(\tx, Ce^{-\varepsilon n})$$
and 
 $$f^n(B_{\varepsilon}(\tx, (r_0/C)e^{-\varepsilon n}))\subset
B_{\varepsilon}(x_j, r_0)\,.$$
This implies that the sequence $\{\UUU_n, n\geq 0\}$ is a uniform family of quasipackings
by Jordan domains.\qed

We note also that the quasipackings we have just constructed have
uniformly bounded overlap by axiom [Degree].

\section{Ahlfors regular conformal dimension and multicurves}

Suppose $f: S^2 \to S^2$ is topologically cxc. By Theorem \ref{thm:qpack_exists}, there exists an associated metric such that  the sequence $\{\UUU_n, n \geq 0\}$ is a family of uniform quasipackings by Jordan domains.   As coverings for the definition of combinatorial moduli, we take $\SSS_n = \UUU_n$.  

\begin{prop} 
\label{prop:mc_and_comb_moduli}
Let $Q >1$, and let $\Gamma$ be a multicurve
with $\lambda(f_{\Gamma,Q}) \ge 1$.
Then, for any $n$ large enough, $\mod_Q([\Gamma],\UUU_n)\gtrsim 1$, where $[\Gamma]$ denotes
the family of all curves in $S^2-P_f$ homotopic to a curve in $\Gamma$.
\end{prop}

\pf Without loss of generality we may assume $\Gamma$ is irreducible.   Write $\Gamma = \{\gamma_1, \ldots, \gamma_m\}$ and equip 
$\R^\Gamma$  with the $L^1$-norm $|\cdot|_1$ norm, so that $|\sum_j a_j \gamma_j |_1 = \sum_j |a_j|$.   For each $1 \leq j \leq m$ choose an annulus $B_j$ 
which is a regular neighborhood of $\gamma_j$ and such that $B_i \intersect B_j = \emptyset, i \neq j$.  
Within each $B_j$ choose a smaller such neighborhood $A_j$ so that $\gamma_j \subset A_j \subset \cl{A}_j \subset B_j$ 
and each inclusion is essential.  
By expansion, there exists a level $n_0$ such that the covering $\UUU_{n_0}$ has the following properties:
\be
\item $s \in \UUU_{n_0}, s \intersect A_j \neq \emptyset \implies s \subset B_j$.  
\item $\mod_Q(A_j, \UUU_{n_0}) >0$ for all $1 \leq j \leq m$.  
\eb

For $n\ge 1$, let $\Gamma_n$ denote the finite family of curves $\wtgamma$ in $[\Gamma]$ arising as connected components of curves of the form $f^{-n}(\gamma_j), \gamma_j \in \Gamma$.  Given such a curve $\wtgamma \subset f^{-n}(\gamma_j)$, denote by $A(\wtgamma)$ the unique component of $f^{-n}(A_j)$ containing $\wtgamma$.  Note that for fixed $n$, the resulting collection of annuli $A(\wtgamma), \wtgamma \in \Gamma_n$, are disjoint.

By the monotonicity and additivity of moduli (Proposition \ref{prop:subadditivity}) we have that 
\[ \mod_Q([\Gamma], \UUU_{n_0+n}) \geq \sum_{\wtgamma\in\Gamma_n} \mod_Q(A(\wtgamma),  \UUU_{n_0+n}).\]

Let $v \in \R^\Gamma$ be the vector of combinatorial moduli at level $n_0$ given by 
\[ v=(\mod_Q(A_1, \UUU_{n_0}), \ldots, \mod_Q(A_m, \UUU_{n_0})).\]

Proposition \ref{prop:covering} implies that if $n \geq 1$ and $\wtgamma \in \Gamma_n$ then 
$$\mod_Q(A(\wtgamma),  \UUU_{n_0+n})= \deg(f:\wtgamma \to f(\wtgamma))^{1-Q}\cdot \mod_Q(A(f(\wtgamma)),  \UUU_{n_0+n-1}).$$
By induction and the fact that degrees multiply under compositions of coverings, for each fixed $n$, the $j$-th entry of the vector
$f_{\Gamma,Q}^{n}(v)$ is the sum of the moduli $\{\mod_Q(A(\wt\gamma)),  \UUU_{n_0+n})\}$ over all curves $\wtgamma\in\Gamma_n$
homotopic to $\gamma_j\in\Gamma$.   By the monotonicity and subadditivity of moduli (Proposition \ref{prop:subadditivity}) we conclude that, for any $n\ge 1$,
$$\mod_Q([\Gamma],\UUU_{n_0+n})\ge |f_{\Gamma,Q}^{n}(v)|_1\,.$$

By the Perron-Frobenius theorem, there is a positive vector $w_Q$ for which 
$f_{\Gamma,Q}(w_Q)=\lambda(f_{\Gamma,Q})\cdot w_Q$.  
By scaling, we may assume $w_Q \leq v$.  Since the entries of the matrix for $f_{\Gamma,Q}$ are  non-negative, we have 
\[ f_{\Gamma,Q}^{n}(v) \geq f_{\Gamma,Q}^n(w_Q) = \lambda(f_{\Gamma,Q})^n w_Q \geq w_Q>0\]
and so 
\[ \liminf_{n \to \infty} |f_{\Gamma,Q}^{n}(v)|_1>0 \]
which completes the proof.
\qed

We conclude with the proof of Theorem \ref{thm:ARconfdim}.  
\gap

\pf By Theorem \ref{thm:qpack_exists}, there is an Ahlfors regular metric $d_\varepsilon\in\GGG(f)$ on 
$S^2$ for which the sequence of coverings $\{\UUU_n, n \geq 0\}$ is a uniform family of quasipackings. 
 Let $\Gamma$ be a multicurve and $[\Gamma]$ the family of all curves homotopic to an element of $\Gamma$.  If $\Gamma$ contains no irreducible multicurve, then $Q(\Gamma)=0 \leq \confdim_{AR}(f)$.   Otherwise, by Lemma \ref{lemma:decreasing} and the definition of $Q(\Gamma)$, we have $\lambda(f_{\Gamma,Q(\Gamma)})=1$ for some $Q(\Gamma)\geq 1$.  By Proposition \ref{prop:mc_and_comb_moduli} applied with $Q=Q(\Gamma)$, 
$\mod_{Q(\Gamma)}([\Gamma], \UUU_n) \gtrsim 1$ as $n \to \infty$.  Since curves in $\Gamma$ are non-peripheral,  there is a positive lower bound for the diameter of any curve in the family $[\Gamma]$. 
 Thus, Corollary \ref{cor:dim} implies that 
$\confdim_{AR}(S^2, d_\varepsilon) \geq Q(\Gamma)$ and so $\confdim_{AR}(f) \geq Q(\Gamma)$.
Since $\Gamma$ is an arbitrary multicurve, we conclude 
$$\confdim_{AR}(f) \geq Q(f)\,.$$

\qed

\appendix
\section{Monotonicity of leading eigenvalues}
\label{secn:appendix}

We first recall some facts concerning non-negative square matrices $A$; see \cite{berman:plemmons}.  \bi

\item {\bf Perron-Frobenius theorem, irreducible version.} \cite[Theorem 1.4]{berman:plemmons}.  A $k$-by-$k$ non-negative matrix $A$  is said to be {\em irreducible} if, for any ordered pair $(i,j)$, $1 \leq i, j \leq k$, there is some power $q >0$  for which $(A^q)_{i,j}>0$.    If $A$ is irreducible, then there is a simple eigenvalue $\lambda(A)$ of $A$ which is larger than the norm of any other eigenvalue, and up to scale, there is a unique corresponding eigenvector, all of whose entries are positive.  

\item {\bf Perron-Frobenius theorem, general version.} \cite[Theorem 1.1]{berman:plemmons}.   If $A$ is merely non-negative, then there exists a non-negative eigenvalue $\lambda(A)$ equal to its spectral radius,  and any corresponding eigenvector is also non-negative.  

\item {\bf Monotonicity} \cite[Corollary 2.1.5]{berman:plemmons}.  The function $A \mapsto \lambda(A)$ satisfies 
\begin{equation}
\label{eqn:monotone}
A \geq B \implies \lambda(A) \geq \lambda(B)
\end{equation}

\item  {\bf Irreducible decomposition} \cite[pp. 39-40]{berman:plemmons}.   
Given any non-negative matrix $A$, there is a permutation matrix $P$ such that  
$U=PAP^{-1}$ has block upper triangular form, where the diagonal blocks $D$ of $U$ are square and either irreducible or zero.  For some diagonal block $D$, $\lambda(D)=\lambda(A)$.  
\ib
\gap

\noindent{\bf Definition.}  Let $p \geq 1$ be an integer.  A {\em Levy cycle of length $p$} is a multicurve $\Gamma=\{\gamma_j,\ j\in\Z/p\Z\}$ such that for each 
$j\in\Z/p\Z$, $f^{-1}(\gamma_j)$ contains a preimage $\delta$ which is homotopic to $\gamma_{j+1}$, and such that $\deg(f: \delta \to \gamma_j)=1$.  

\begin{lemma}
\label{lemma:no_levy} 
If $f: S^2 \to S^2$ is a branched covering satisfying Axiom [Expansion] with respect to an open covering $\UUU_0$, then $f$ has no Levy cycles.  
\end{lemma}

\pf  Fix a metric on the sphere compatible with its topology.  
Axiom [Expansion] implies that there are constants $d_n \downarrow 0$ as $n \to \infty$ such that 
$\max\{\diam U | U \in \UUU_n\} \leq d_n$.

Suppose $f$ had a Levy cycle $\Gamma$ of length $p$, and let $\gamma \in \Gamma$.  
Then $g=f^p$ also satisfies Axiom [Expansion] with respect to $\UUU_0$, and $g^{-1}(\gamma)$ 
has a connected component $\delta$ homotopic to $\gamma$ and satisfying $\deg(g:\delta \to \gamma)=1$.  
There is an open annulus $A \subset S^2-P_f$ containing $\gamma$ such that 
the inclusion map $\gamma \hookrightarrow A$ is essential.  
By compactness and Axiom [Expansion], there exists $n_0 \in \N$ such that 
\[ \UUU_{n_0}(\gamma) = \{U \in \UUU_{n_0} | U \intersect \gamma \neq \emptyset\} \subset A.\]
Let $N=\#\UUU_{n_0}(\gamma)$.  Since $\deg(g:\delta \to \gamma)=1$ and $\delta$ is homotopic to $\gamma$, it follows by induction and the construction of the annulus $A$ that for all $k \in \N$, there exists a  component $A_k$ of $g^{-k}(A)$ homotopic to $A$ such that $\deg(g^k: A_k \to A)=1$, i.e. $g^k|A_k$ is a homeomorphism onto $A$.  Hence, the annulus $A_k$ contains a unique component $\delta_k$ of $g^{-k}(\gamma)$ which is homotopic to $\gamma$, and $\delta_k$ is covered by $N$ elements of $\UUU_{n_0+pk}$.  Thus $\diam \delta_k \leq N\cdot d_{n_0+pk} \to  0$ as $k \to \infty$.  But this is impossible:  
since $\gamma$ is non-peripheral, there is a positive lower bound on the diameter of any curve homotopic to $\gamma$.

\qed

\begin{lemma}
\label{lemma:decreasing}
Let  $f: S^2 \to S^2$ be a branched covering satisfying  Axiom [Expansion], and let $\Gamma$ be a multicurve.
\be
\item If $\Gamma$ does not contain an irreducible multicurve, then $f_{\Gamma,Q}$ is nilpotent and $\lambda(f_{\Gamma,Q})=0$ for all $Q\ge 1$.
\item If $\Gamma$ contains an irreducible multicurve, then $\lambda(f_{\Gamma,1}) \geq 1$, and the function $Q \mapsto \lambda(f_{\Gamma,Q})$ is strictly decreasing on $[1,\infty)$ and tends to zero as $Q$ tends to $\infty$.
\eb
\end{lemma}

\pf  Re-indexing the elements of $\Gamma$, we may assume the matrix $(f_{\Gamma,Q})$ has block-upper 
triangular form for all $Q\ge 1$. 

Suppose that  $\Gamma$ contains no irreducible multicurve. 
Then 
the matrix $(f_{\Gamma,Q})$ is upper triangular and has zeros on the diagonal.  
Hence $f_{\Gamma,Q}$ is nilpotent and 
$\lambda(f_{\Gamma,Q})=0$.

Suppose now that $\Gamma$ contains an irreducible multicurve $\Gamma'$.
It follows that
$\lambda(f_{\Gamma,1}) \geq \lambda(f_{\Gamma',1})$.
Since the matrix $(f_{\Gamma', 1})$ is irreducible, non-negative and with positive
entries at least $1$, there is some permutation matrix $P'$ and a re-indexing such that
$$f_{\Gamma',1} \ge \mtwo{P'}{0}{0}{0}\,.$$
We then have 
 by Equation (\ref{eqn:monotone}) 
\[ \lambda(f_{\Gamma,1}) \geq \lambda(f_{\Gamma',1}) \geq \lambda(P')=1.\]

We now prove the second assertion.   For convenience, denote by $A_Q$ the matrix $(f_{\Gamma,Q})$.  If $Q_1>Q_2$ then $A_{Q_1} \leq A_{Q_2}$ entrywise and  Equation (\ref{eqn:monotone}) 
implies that $\lambda(A_{Q_1}) \leq \lambda(A_{Q_2})$.    
If equality holds for distinct $Q_1, Q_2$, then $\lambda(A_Q)$ is constant for all $Q_2 \leq Q \leq Q_1$.  
Since eigenvalues are algebraic functions,  this would imply $\lambda(A_Q)$ is constant for all $Q$.  
Since $\lambda(A_1)\geq 1$ by assumption, it suffices to show that $\lambda(A_Q) \to 0$ as $Q \to \infty$.  

The definition of $f_{\Gamma,Q}$ implies that 
\[ A_Q = B_Q + C\]
where $B_Q, C$  are  non-negative, $B_Q \to 0$, and the entries of $C$ are of the form $1^{1-Q}+\ldots + 1^{1-Q}$.  
Hence $C$ is constant in $Q$.  

In this paragraph, we prove that $C^m=0$ where $m=\#\Gamma$.  Suppose  $D$ is an irreducible diagonal block in the decomposition of $C$.  Then $(D^q)_{i,i}>0$ for some index $i$ and some power $q>0$.  But this implies that $\Gamma$ contains a Levy cycle, which is impossible by Lemma \ref{lemma:no_levy}.    
Hence all diagonal blocks are zero, which implies $C^m=0$.  

Since $C^m=0$, every term in the expansion of $(B_Q+C)^m$ contains $B_Q$ as a factor.  Therefore 
\[ \lim_{Q \to \infty} A^m_Q =  \lim_{Q \to \infty} (B_Q+C)^m = 0\]
entrywise.   Hence $\lambda(A_Q^m) = \lambda(A_Q)^m \to 0$ and so $\lambda(A_Q) \to 0$.  
\qed

\def\cprime{$'$}

\end{document}